\documentclass[11pt]{article}

\usepackage{amsmath, amssymb}
\usepackage{latexsym}
\usepackage{amsfonts}
\usepackage{graphicx}
\usepackage{amscd}
\usepackage{amsbsy}
\newtheorem{thm}{Theorem}[section]
\newtheorem{cor}[thm]{Corollary}
\newtheorem{lem}[thm]{Lemma}

\newcommand{\be}{\begin{equation}}
\newcommand{\ee}{\end{equation}}

\newcommand{\pf}{\medskip \noindent {\sl Proof}. \ }
\newcommand{\qed}{\hfill $\Box$ \\ \medskip}
\def\nm{\noalign{\medskip}}
\newcommand{\Om}{\Omega}

\newcommand{\la}{\langle}
\newcommand{\ra}{\rangle}

\newcommand{\ds}{\displaystyle}
\newcommand{\p}{\partial}

\newcommand{\pd}[2]{\frac {\p #1}{\p #2}}

\newcommand{\eqnref}[1]{(\ref {#1})}
\newcommand{\na}{\nabla}
\newcommand{\ep}{\epsilon}
\newcommand{\vp}{\varphi}
\newcommand{\Scal}{\mathcal{S}}
\newcommand{\Dcal}{\mathcal{D}}
\newcommand{\Kcal}{\mathcal{K}}

\newcommand{\RR}{\mathbb{R}}

\newcommand{\NN}{\mathbb{N}}

\begin{document}
\title{Optimal Estimates for the Electric Field in Two-Dimensions}

\author{Habib Ammari (ammari@cmapx.polytechnique.fr)\thanks{Centre
de Math{\'e}matiques Appliqu{\'e}es, CNRS UMR 7641 and Ecole
Polytechnique, 91128 Palaiseau Cedex, France.}\\
Hyeonbae Kang (hkang@math.snu.ac.kr)\thanks{School of Mathematical
Sciences, Seoul National University, Seoul 151-747, Korea.} \\
Hyundae Lee (hdlee@math.snu.ac.kr)\footnotemark[2] \\
Jungwook Lee (wook9@snu.ac.kr)\footnotemark[2] \\
Mikyoung Lim (mklim@cmapx.polytechnique.fr)\footnotemark[1] }

\date{}
\maketitle

\begin{abstract}
We establish both upper and lower bounds on the electric field in
the case where two circular conductivity inclusions are very close
but not touching. We also obtain such bounds when a circular
inclusion is very close to the boundary of a circular domain which
contains the inclusion. The novelty of these estimates, which
improve and make complete our earlier results in \cite{AKL}, is
that they give an optimal information about the blow-up of the
electric field as the conductivities of the inclusions degenerate.

\vskip 0.5\baselineskip

\end{abstract}

\noindent {\footnotesize Mathematics subject classification
(MSC2000): Primary  35J25; Secondary 73C40}



\section{Introduction and statements of results}
The purpose of this paper is to set out optimal gradient estimates
for solutions to the isotropic conductivity problem in the
presence of adjacent conductivity inclusions as the distance
between the inclusions goes to zero and their conductivities
degenerate. This difficult question arises in the study of
composite media. Frequently in composites, the inclusions are very
closely spaced and may even touch; see \cite{bab}. It is quite
important from a practical point of view to know whether the
electric field (the gradient of the potential) can be arbitrarily
large as the inclusions get closer to each other or to the
boundary of the background medium.

There have been some important works on the estimates of the
gradient of the solution to the conductivity problem in the
presence of inclusions. For finite and strictly positive
conductivities, it was shown by Bonnetier and Vogelius in
\cite{bonvog} that the gradient of $u$ remains bounded for
circular touching inclusions of comparable radii. Li and Vogelius
showed in \cite{livog} that $\nabla u$ is bounded independently of
the distance between the inclusions $B_1$ and $B_2$, provided that
the conductivities stay away from $0$ and $+ \infty$.  It is worth
mentioning that the result of \cite{livog} is much more general:
it holds for arbitrary number of inclusions with arbitrary shape.
This result has been recently extended to elliptic systems by Li
and Nirenberg in \cite{linir}. On the other hand, for two
identical perfectly conducting circular inclusions (with
$k_1=k_2=+\infty$) which are $\epsilon$ apart, it has been shown
in \cite{budcar} (see also \cite{mar} and \cite{keller}) that the
gradient generally becomes unbounded as the distance $\epsilon$
approaches zero. The rate at which this gradient becomes unbounded
has actually been calculated in \cite{budcar}, for a special
solution. For this special solution, the rate turns out to be
$\epsilon^{-1/2}$. In \cite{AKL}, a lower bound for the gradient
of the solution to the conductivity problem for arbitrary
conductivities, possibly degenerating, has been obtained. One of
our objectives in this paper is to improve this bound and complete
it by deriving an optimal upper one as well.

In this paper we consider the following two situations: when two
circular conductivity inclusions are very close but not touching
and when a circular inclusion is very close to the boundary of the
domain where the inclusion is contained. These two simple
two-dimensional models illustrate very well the feature of our
estimates. We believe that they extend to arbitrary-shaped
inclusions if their contact reduces to a point.

To describe the first situation we consider in this paper, let
$B_1$ and $B_2$ be two circular inclusions contained in a matrix
which we assume to be the free space $\RR^2$. For $i=1,2$, we
suppose that the conductivity $k_i$ of the inclusion $B_i$ is a
constant different from the constant conductivity of the matrix,
which is assumed to be $1$ for convenience. The conductivity $k_i$
of the inclusion may be $0$ or $+\infty$. The zero conductivity
indicates that the inclusion is an insulated inclusion while the
infinite conductivity indicates a perfect conductor. We are
especially interested in the case of extreme conductivities $k_i
\rightarrow +\infty$ or $k_i \rightarrow 0$.

Given an entire harmonic function $H$, the first conductivity
problem we consider in this paper is the following:
 \be \label{cond-trans}
 \begin{cases}
 \ds \nabla  \cdot \Bigr( 1 + \sum_{i=1,2} (k_i -1) \chi(B_i)
 \Bigr) \nabla u =0 \quad \mbox{in } \RR^2, \\
 u(X)- H(X) = O(|X|^{-1}) \quad \mbox{as } |X| \to + \infty.
 \end{cases}
 \ee
The electric field is given by $\nabla u$, where $u$ is the
solution to \eqnref{cond-trans} and represents the perturbation of
the field $\nabla H$ in the presence of the two inclusions $B_1$
and $B_2$. For applications to the theory of composite materials,
it is particularly important to consider the case when $\nabla H$
is a uniform field, {\em i.e.}, $H(X)=A \cdot X$ for some constant
vector $A$. The equation \eqnref{cond-trans} can be rewritten in
the following form to emphasize the transmission conditions on $\p
B_i$, $i=1,2$:
$$
\begin{cases}
\ds \Delta u =0  & \mbox{in } \Om \setminus (\p B_1 \cup \p
B_2), \\
\ds u|_+ = u|_-  & \mbox{on } \p B_i, \ i=1,2, \\
\nm
\ds \pd{u}{\nu} \bigg |_{+} = k_i \pd{u}{\nu} \bigg |_{-} & \mbox{on } \p
B_i, \
i=1,2,  \\
u(X)- H(X) = O(|X|^{-1}) \quad &\mbox{as } |X| \to + \infty.
\end{cases}
$$
Here and throughout this paper the subscript $\pm$ indicates the
limit from outside and inside the domain, respectively. If
$k_i=0$, then the transmission condition on the normal derivatives
of $u$ should be replaced with $\pd{u}{\nu} |_{+}=0$ on $\p B_i$,
while if $k_i= + \infty$, it should be replaced with
$u=\mbox{constant}$ on $B_i$.

As has already been said, we are interested in the behavior of the
gradient of the solution to the equation \eqnref{cond-trans} as
the distance between the inclusions $B_1$ and $B_2$ goes to zero
for arbitrary conductivities $k_1$ and $k_2$, possibly
degenerating.

To state the first main result of this paper, let us first fix
notation. For $i=1,2$, let $B_i=B(Z_i, r_i)$, the disk centered at
$Z_i$ and of radius $r_i$. Let $R_i$, $i=1,2$, be the reflection
with respect to $\p B_i$, {\em i.e.},
$$
R_i (X): = \frac{r_i^2 (X-Z_i)}{|X-Z_i|^2}+Z_i, \quad i=1,2.
$$
It is easy to see that the combined reflection $R_1R_2$ and
$R_2R_1$ have unique fixed points. Let $I$ be the line segment
between these two fixed points. Let $X_j$, $j=1,2$, be the point
on $\p B_j$ closest to the other disk. We also let
$$
r_{\min} := \min(r_1,r_2), \quad r_{\max}:=\max(r_1,r_2), \quad
r_* :=\sqrt{(2r_1r_2)/(r_1+r_2)},
$$
$$
\lambda_i := \frac{k_i+1}{2(k_i -1)}, \ \quad i=1,2 \quad
\mbox{and} \quad \tau:= \frac{1}{4\lambda_1 \lambda_2}.
$$
We obtain the following result on the blow-up of the gradient.
\begin{thm} \label{firstmain}
Let $\epsilon := \mbox{dist}(B_1,B_2)$ and let $\nu^{(j)}$ and
$T^{(j)}$, $j=1,2$, be the unit normal and tangential vector
fields to $\p B_j$, respectively. Let $u$ be the solution of
\eqnref{cond-trans}.
\begin{itemize}
\item[{\rm (i)}] If $\ep$ is sufficiently small, there is a constant $C_1$
independent of $k_1$, $k_2$,
$r_1$, $r_2$, and $\ep$ such that
\begin{equation} \label{estthm}
\frac{C_1 \ds \inf_{X \in I} |\la \nabla H(X), \nu^{(j)}(X_j) \ra
|}{1 -\tau +    (r_*/r_{\min}) \sqrt{\epsilon}} \leq |\nabla u|_{+}
(X_j) |, \quad j=1,2,
\end{equation}
provided that $k_1, ~k_2>1$, and
\begin{equation} \label{estthm1}
\frac{C_1 \ds \inf_{X \in I} |\la \nabla H(X), T^{(j)}(X_j) \ra
|}{1 -\tau +    (r_*/r_{\min}) \sqrt{\epsilon}} \leq |\nabla u
|_{+} (X_j) |, \quad j=1,2,
\end{equation} provided that $k_1, ~k_2<1$.
\item[{\rm (ii)}] Let $\Om$ be
a bounded set containing $B_1$ and $B_2$. Then there is a constant $C_2$
independent of $k_1$, $k_2$,
$r_1$, $r_2$, $\ep$, and $\Om$ such that
\begin{equation} \label{estthm-upper}
\|\nabla u \|_{L^\infty(\Om)} \leq
\frac{C_2 \| \nabla H \|_{L^\infty(\Om)}}{1 -|\tau| +  ( r_*/r_{\max})
\sqrt{\epsilon}}.
\end{equation}
\end{itemize}
\end{thm}

Note that if $H(X)=A\cdot X$ for some constant vector $A$, which
is the most interesting case, then the quantity
$$
\la \nabla H(X), \nu^{(j)}(X_j) \ra = \la A,\nu^{(j)}(X_j) \ra,
$$
and hence it does not vanish if we choose $A$ appropriately.

Theorem \ref{firstmain} quantifies the behavior of $\nabla u$ in
terms of the conductivities of the inclusions, their radii, and
the distance between them. For example, if $k_1$ and $k_2$
degenerate to $+\infty$ or zero, then $\tau=1$ and hence
\eqnref{estthm} and \eqnref{estthm-upper} read
\begin{equation} \label{estthm-1}
\frac{C^\prime_1}{(r_*/r_{\min}) \sqrt{\epsilon}} \leq |\nabla
u(X_j) |, \quad j=1,2,
  \quad \|\nabla u \|_{L^\infty(\Om)} \leq  \frac{C^\prime_2}{( r_*/r_{\max})
\sqrt{\epsilon}},
\end{equation}
for some positive constants $C^\prime_1$ and $C^\prime_2$, which
shows that $\nabla u$ blows up at the rate of $\epsilon^{-1/2}$ as
the inclusions get closer. It further shows that the gradient
blows up at $X_1$ and $X_2$, $X_j$ for $j=1,2$, being the point on
$\p B_j$ closest to the other disk.

The lower bounds in \eqnref{estthm} and \eqnref{estthm1} are
improved versions of the one obtained in \cite{AKL}. The proofs of
our new estimates make use of quite explicit but nontrivial
expansion formulae, originally derived in \cite{AKKL}. They are
achieved by using a significantly different method from
\cite{bonvog}, \cite{livog}, and \cite{budcar}.

Another interesting situation is when the inclusion is very close
to the boundary. To describe this second situation, suppose that
$\Omega$, which is a disk of radius $\rho$, contains an inclusion
$B$, which is a disk of radius $r$. Suppose also that the
conductivity of $\Om$ is $1$ and that of $B$ is $k \neq 1$. The
conductivity problems considered in this case are the following
Dirichlet and Neumann problems: for a given $f \in
\mathcal{C}^{1,\alpha}(\p \Om)$, $\alpha >0$,
\begin{equation}\label{eq-diri}
\left\{ \begin{array}{ll}
\nabla \cdot (1+(k-1)\chi(B))\nabla u = 0 &~~ \mbox{in}~~ \Omega, \\
  u=f &~~ \mbox{on}~~ \partial\Omega,
\end{array}\right.
\end{equation}
and for a given $g \in \mathcal{C}^\alpha(\p \Om)$
\begin{equation}\label{eq-neu}
\left\{ \begin{array}{ll}
\nabla \cdot (1+(k-1)\chi(B))\nabla u = 0 &~~ \mbox{in}~~ \Omega, \\
\nm
\ds \pd{u}{\nu} =g &~~ \mbox{on}~~ \partial\Omega~.
\end{array}\right.
\end{equation}
To ensure existence and uniqueness of a solution to
\eqnref{eq-neu}, we suppose that $\int_{\p\Om} g=0$ and
$\int_{\p\Om} u=0$. Let $X_1$ be the point on $\p B$ closest to
$\p\Om$ and $X_2$ be the point on $\p \Om$ closest to $\p B$, and
let $R_B$ and $R_\Om$ are reflections with respect to $\p B$ and
$\p \Om$, respectively. Let $P_1$ and $P_2$ be fixed points of
$R_BR_\Om$ and $R_\Om R_B$, respectively, and let $J_1$ be the
line segment between $P_1$ and $X_1$ and $J_2$ that between $P_2$
and $X_2$.

The second main result of this paper is the following triplet of
estimates for the gradient of the solutions to \eqnref{eq-diri}
and \eqnref{eq-neu}. Let $\Dcal_\Om(f)$ and $\Scal_\Om(g)$ denote
the double and single layer potentials whose definitions are given
in Section 2. For the Dirichlet problem we have the following
theorem.

\begin{thm}\label{secondmain}
Let
$$
\epsilon:= \mbox{dist}(B, \p\Om), \quad \sigma:= \frac{k-1}{k+1},
\quad
  r^*:=\sqrt{\frac{\rho-r}{\rho r}},
$$
and let $u$ be the solution to \eqnref{eq-diri}.
\begin{itemize} \item[{\em (i)}]
If $k>1$, then there exists a constant $C_1$ independent of $k$,
$r$, $\epsilon$, and $f$ such that for $\epsilon$ small enough,
\begin{equation}\label{est}
\frac{C_1 \ds\inf_{X \in J_1} | \la \nabla \Dcal_\Om(f)(X), \nu_B
(X_1) \ra |}{1-\sigma
   + 4 r^*\sqrt{\epsilon}}\leq |\na
u|_{+}(X_1) |,
\end{equation}
and
\begin{equation}\label{est-100}
\frac{C_1 \ds\inf_{X \in J_2} | \la \nabla \Dcal_\Om(f)(X),
\nu_\Omega (X_2) \ra |}{1-\sigma
   + 4 r^*\sqrt{\epsilon}}\leq |\na
u|_{-}(X_2) |.
\end{equation}
Here $\nu_B$ and $\nu_\Omega$ denote the outward unit normal to
$\p B$ and $\p\Om$.
\item[{\rm (ii)}] For any $k \neq 1$, there exists a constant
$C_2$ independent of $k$, $r$, and $\epsilon$ such that for
$\epsilon$ small enough,
\begin{equation}\label{est1}
\|\na u\|_{L^\infty(\Om)}\leq \ds\frac{C_2 \| f
\|_{\mathcal{C}^{1,\alpha}(\p \Om)}}{1- |\sigma| +
r^*\sqrt{\epsilon}}.
  \end{equation}
\end{itemize}
\end{thm}

For the Neumann problem the following theorem holds.
\begin{thm}\label{thirdmain}
Let $\epsilon$, $\sigma$, $r^*$ be defined as in Theorem
\ref{secondmain}.
\begin{itemize}
\item[{\rm (i)}] If $k<1$, then there exists a constant $C_1$
independent of $k$, $r$, $\epsilon$, and $g$ such that for
$\epsilon$ small enough,
\begin{equation}\label{est-neu}
\frac{C_1 \ds\inf_{X \in J_1} | \la \nabla \Scal_\Om(g)(X),
T_B(X_1) \ra |}{1+\sigma
   + 4 r^*\sqrt{\epsilon}}\leq |\na
u|_{+}(X_1) |,
\end{equation}
and
\begin{equation}\label{est-neu-100}
\frac{C_1 \ds\inf_{X \in J_2} | \la \nabla \Scal_\Om(g)(X),
T_\Omega (X_2) \ra |}{1+\sigma
   + 4 r^*\sqrt{\epsilon}}\leq |\na
u|_{-}(X_2) | .
\end{equation}
Here $T_B$ and $T_\Omega$ denote the positively oriented unit
tangent vector fields on $\p B$ and $\p\Om$, respectively.

\item[{\rm (ii)}] For any $k \neq 1$, there exists a constant
$C_2$ independent of $k$, $r$, and $\epsilon$ such that for
$\epsilon$ small enough,
\begin{equation}\label{est1-neu}
\|\na u\|_{L^\infty(\Om)}\leq \ds\frac{C_2 \| g
\|_{\mathcal{C}^{\alpha}(\p \Om)}}{1-|\sigma| +
r^*\sqrt{\epsilon}}. \ee
\end{itemize}

\end{thm}

If $Z$ is the center of $\Om$ and if $f(X)= A\cdot X$ for some
constant vector $A$, then $\Dcal_\Om (f)(X) = \frac{1}{2} A\cdot
X$ for $X \in \Om$ and $\Dcal_\Om (f)(X) = -\frac{\rho A \cdot
X}{2|X-Z|^2}$ for $X \in \RR^2 \setminus \overline{\Om}$, and
hence we can achieve $$ \la \nabla \Dcal_\Om(f)(X), \nu_B(X_1) \ra
\neq 0 \mbox{ and } \la \nabla \Dcal_\Om(f)(X), \nu_\Omega(X_2)
\ra \neq 0 \quad \mbox{for any } X, $$ by choosing $A$
appropriately. Likewise, if $g := A \cdot \nu$ on $\p\Om$, then
$\Scal_\Om (g)= -\frac{1}{2} A\cdot X + \mbox{constant}$.

Theorem \ref{secondmain} shows that in the case of the Dirichlet
problem, if the inclusion is a perfect conductor ($k =+ \infty$
and hence $\sigma=1$), then
$$
  \frac{C^\prime_1}{r^*\sqrt{\epsilon}}\leq  \| \na u \|_{L^\infty(\Om)}\leq
\ds\frac{C^\prime_2}{r^*\sqrt{\epsilon}},
$$
for some positive constants $C^\prime_1$ and $C^\prime_2$. Thus
$\nabla u$ blows up at the rate of $\ep^{-1/2}$ as long as the
magnitude of $r$ is much larger than that of $\ep$. It also shows
that the gradient blows up at the points $X_1$ and $X_2$. On the
other hand, for the Neumann problem, according to Theorem
\ref{thirdmain} the situation is reversed: $\nabla u$ blows up for
an insulator. If $r$ is of the same order as $\ep$, then $r^*
\approx \frac{1}{\sqrt{\ep}}$ and hence $\nabla u$ does not blow
up. In fact, it stays bounded and an asymptotic expansion of the
solution as $\ep \to 0$ can be derived. See for instance
\cite{AAK, AK2,AKT} for this.

One may think that Theorems \ref{secondmain} and \ref{thirdmain}
can be easily derived from Theorem \ref{firstmain} by reflection
or conformal mapping. This is, as it will be shown, far from being
true. The proofs of Theorems \ref{secondmain} and \ref{thirdmain}
require delicate analysis and careful and tricky estimates.

In this paper we only deal with the two dimensional case. It seems
challenging to obtain similar results in three dimensions. At this
moment it is even not clear what the blow-up rate of the gradient
would be in three dimensions. In this direction, we have recently
established in \cite{ADKL} that, unlike the two-dimensional case,
if the inclusions are grounded conductors (the Dirichlet boundary
condition on the inclusions is set to be zero) and of spherical
shape, the gradient stays, to our surprise, bounded regardless of
the separation distance between them.

This paper is organized as follows. In Section 2 we review some
preliminary facts on layer potentials and representations of the
solution to the conductivity problem obtained in \cite{AKL}.
Theorem \ref{firstmain} is proved in Section 3, and Theorems
\ref{secondmain} and \ref{thirdmain} in Section 4. Although our
results hold for special cases, we believe that they extend to
arbitrary-shaped inclusions if their contact reduces to a point.

\section{Preliminaries}

To make our paper self-contained and our exposition clear, we
introduce our main tools for studying the conductivity problems
and collect some preliminary results regarding layer potentials.
We also linger over a description of some results from our earlier
papers. The material in Lemmas \ref{represen} and \ref{conj-2} is
however, up to our knowledge,  new.

Let $D$ be a bounded Lipschitz domain in $\RR^2$ and $\Scal_D
\phi$ and $\Dcal_D \phi $ denote the single and double layer
potentials of a function $\phi \in L^2(\partial D)$, namely,
\begin{align*}
\Scal_D \phi(X) & = \frac{1}{2\pi} \int_{\p D} \ln
|X-Y| \phi(Y) \, d\sigma(Y), \quad X \in \RR^2,\\
\Dcal_D \phi(X) & = \frac{1}{2\pi} \int_{\p D} \frac{\la
Y-X,\nu_Y\ra}{|X-Y|^2} \phi(Y) \, d\sigma(Y), \quad X \in
\RR^2\setminus\p D.
\end{align*}
 For a function $u$ defined on $\RR^2 \setminus \p D$, we set
$$ \pd{u}{\nu} \bigg |_\pm(X) := \lim_{t \to 0^+} \la \nabla u(X
\pm t \nu_X), \nu_X \ra\;, \quad X \in \p D \;, $$ if the limits
exist. Here and throughout this paper $\nu_X$ is the outward
normal to $\p D$ at $X$ and $\la~ , \ra$ denotes the scalar
product in $\RR^2$. The following jump relations for the single
and double layer potentials are well-known \cite{folland76,
verchota84}.

\begin{lem}\label{jump}
For $\phi\in L^2(\p D)$ we have
\begin{align}
\pd{}{\nu}\Scal_D\phi\Big|_\pm(X)
&=\big(\pm\frac{1}{2}+\Kcal^*_D\big)\phi(X) \quad \mbox{a.e.}~
X\in\p D, \label{singlejump} \\
(\Dcal_D\phi)|_
\pm(X)&=(\mp\frac{1}{2}+\Kcal_D\big)\phi(X) \quad \mbox{a.e.}~ X\in\p
D, \label{doublejump}
\end{align}
where $\Kcal_D$ is defined by
$$
\Kcal_D\phi(X)=\frac{1}{2\pi} \mbox{p.v.} \int_{\p D}\frac{\la
Y-X,\nu_Y\ra}{|X-Y|^2}\phi(Y)d\sigma(Y),
$$
and $\Kcal^*_D$ is the $L^2$-adjoint of  $\Kcal_D$, i.e.,
$$
\Kcal^*_D\phi(X)=\frac{1}{2\pi} \mbox{p.v.} \int_{\p D} \frac{\la
X-Y,\nu_X\ra}{|X-Y|^2}\phi(Y)d\sigma(Y).
$$
Here p.v. denotes the Cauchy principal value.
\end{lem}

Note that if $D$ is a two-dimensional disk with radius $r$, then
$$
\frac{\la X-Y,\nu_X\ra}{|X-Y|^2}=\frac{1}{2r}, \quad \forall X,Y\in \p
D,~X\ne Y,
$$
and hence
\be \label{KDKD}
\Kcal^*_D\phi(X)=\Kcal_D\phi(X)=\frac{1}{4\pi r}\int_{\p D}\phi(Y)d\sigma,
\quad X \in \p D.
\ee

When $D$ is a disk, we define the function $R_D(f)$ for a function
$f$ by
$$
R_D(f)(X): =f(R_D(X)).
$$
Using the jump formula for the single layer potential, we
get the following lemma \cite{AKKL}.
\begin{lem}\label{lm1}
Let $D$ be a disk in $\RR^2$ and let $R_D$ denote the reflection
with respect to $\p D$. If $v$ is a harmonic in $D$ and continuous
on $\overline{D}$ then
\begin{equation}\label{lm12}
\Scal_{D}\left(\pd{v}{\nu} \bigg |_{\p D}
\right)(X)= -\frac{1}{2}(R_{D} v)(X)+C , \quad
X\in\RR^2\setminus\overline{D} ,
\end{equation}
where $C$ is some constant. Analogously, if $v$ is harmonic in
$\RR^2\setminus\overline{D}$, continuous on $\RR^2\setminus{D}$,
and $v(X)\rightarrow~0$  as $  |X|\rightarrow 0$, then
\begin{equation}\label{lm13}
\Scal_{D}\left( \pd{v}{\nu} \bigg |_{\p D} \right)(X)= \frac{1}{2}
(R_{D} v)(X)+C , \quad X\in\overline{D},
\end{equation}
for some constant $C$.
\end{lem}

We need the following lemma which was first proved in \cite{AKL}.
\begin{lem}\label{AKL}
Suppose that $B_1$ and $B_2$ are two disjoint disks and let
$R_i=R_{B_i}$ be the reflection with respect to $\p B_i$. Then the
solution to \eqnref{cond-trans} is represented as
\begin{equation} \label{repp}
u(X) = H(X) + \Scal_{B_1} \vp_1(X) +  \Scal_{B_2} \vp_2 (X),\quad
X \in\Om,
\end{equation}
where $\vp_i \in L^2_0(\p B_i)$, $i=1,2$, is the unique solution
to the system of integral equations \be \label{multi-int-2}
\lambda_l \varphi_l- \pd{(\Scal_{B_i}\varphi_i)}{\nu^{(l)}} \big
|_{\p B_i}=\pd{H}{\nu} \big |_{\p B_i} \quad \mbox{on } \p B_l,
\quad l=1, 2, \ i \neq l, \ee with $\lambda_i =
\frac{k_i+1}{2(k_i-1)}$ and $\nu^{(l)}$ is the outward unit normal
to $\p B_l$. Moreover, the potentials $\vp_1$ and $\vp_2$ are
explicitly given by
\begin{equation} \label{phi_blow_up_1}
\begin{array}{rcl}
\ds \varphi_1 &=& \ds \frac{1}{\lambda_1} \sum_{m=0}^{+\infty}
\frac{1}{(4\lambda_1 \lambda_2)^{m}} \pd{}{\nu^{(1)}}
               \left [(R_2 R_1)^m (I-\frac{1}{2\lambda_2} R_2 )H \right ]
\Bigr |_{\p B_1}, \\
\nm \ds \varphi_2 &=& \ds \frac{1}{\lambda_2} \sum_{m=0}^{+\infty}
\frac{1}{(4\lambda_1 \lambda_2)^{m}} \pd{}{\nu^{(2)}}
               \left [(R_1 R_2)^m (I - \frac{1}{2\lambda_1}R_1)H\right ]
\Bigr |_{\p B_2},
\end{array}
\end{equation}
where the series in \eqnref{phi_blow_up_1} converge absolutely and
uniformly.
\end{lem}

In Lemma \ref{AKL}, the space  $L^2_0(\p B_i)$ denotes the set of
all $g \in L^2(\p B_i)$ having mean value zero: $\int_{\p B_i}
g=0$. The following lemma from \cite{AKL} is also of use to us.

\begin{lem}\label{conj}
Let $u$ be the solution of \eqnref{cond-trans} and let $\widetilde H$
be a harmonic conjugate to $H$. Let $v$ be the solution to the
conductivity problem
\begin{equation} \label{eqv}
\ \left \{
\begin{array}{l}
\ds \nabla  \cdot \Bigr( 1 + \sum_{i=1,2} (\frac{1}{k_i} -1)
\chi(B_i) \Bigr) \nabla v =0 \quad \mbox{in } \RR^2, \\ \nm \ds
v(X)- \widetilde H(X) = O(|X|^{-1}).
\end{array}
\right .
\end{equation}
Then
$$
\frac{\partial u}{\partial  T} = - \frac{\partial v}{
\partial \nu^{(i)}} \bigg |_+  \quad \mbox{ on }~ \p B_i ,~
i=1,2.
$$
\end{lem}

We now turn our attention to the second situation, {\it i.e.}
problems \eqnref{eq-diri} and \eqnref{eq-neu}, when both $\Om$ and
$B$ are disks. We first note that $\Dcal_\Om f$ is
$\mathcal{C}^{1,\alpha}$ in $\overline{\Om}$ and $\RR^2 \setminus
\Om$ since $f \in \mathcal{C}^{1,\alpha}(\p\Om)$. It was shown in
\cite{KS96, KS2000} that the solution $u$ to the problem
\eqnref{eq-diri} for a fixed Dirichlet data $f$ is given by
\begin{equation}\label{cd}
  u(X)=\Dcal_\Omega(f)(X)-\Scal_\Omega(g)(X)+\Scal_B(\vp)(X), \quad X \in
\Om , \quad g:=\frac{\p u}{\p\nu}
  \big |_{\p\Om},
\end{equation}
where $\vp$ with mean value zero satisfies the integral equation
$$ (\lambda I- \Kcal_B^*) \vp = \pd{}{\nu} \big
(\Dcal_\Omega(f)-\Scal_\Omega(g) \big ) \quad \mbox{on } \p B, $$
with $\lambda=\frac{\ds k+1}{\ds 2(k-1)}$. Since $B$ is a disk, it
follows from \eqnref{KDKD} that $\Kcal_B^* \vp \equiv 0$ on
$L^2_0(\p B)$ and hence
 \be \label{vpform}
 \lambda \vp = \pd{}{\nu} \big
 (\Dcal_\Omega(f)-\Scal_\Omega(g) \big ) \quad \mbox{on } \p B.
 \ee
On the other hand, $g=\frac{\p u}{\p\nu}|_{\p\Om}$ yields
$$
g = \pd{}{\nu} \big
(\Dcal_\Omega(f)-\Scal_\Omega(g)+\Scal_B(\vp) \big ) \big |_{-} \quad
\mbox{on } \p \Om.
$$
Since $\pd{}{\nu} \Scal_\Omega(g)|_-= (-\frac{1}{2} I +
\Kcal_\Om^*) g$ and $\Om$ is a disk, $\pd{}{\nu}
\Scal_\Omega(g)|_-= -\frac{1}{2} g$ on $\p B$. Thus we get
\be \label{gform}
\frac{1}{2} g = \pd{}{\nu} \big
(\Dcal_\Omega(f)+\Scal_B(\vp) \big ) \big |_{-} \quad
\mbox{on } \p \Om.
\ee
It then follows from \eqnref{vpform} and \eqnref{gform} that $g$
and $\vp$ are the solution of the following system of integral
equations
\begin{equation}\label{sy}
\left\{
\begin{array}{rcl}
\ds\frac{1}{2}g - \frac{\p (\Scal_B \vp)}{\ds \p\nu_\Om} & = &\ds
\frac{\p(\Dcal_\Om f)}{\p\nu_\Om}
\quad \mbox{on } \p \Om, \\
\nm
\ds \lambda\vp + \frac{\p (\Scal_\Om g)}{\p\nu_B}& =& \ds \frac{\p
(\Dcal_\Om
f)}{\p\nu_B} \quad \mbox{on } \p B .
\end{array} \right.
\end{equation}
Observe the similarity of \eqnref{sy} to \eqnref{multi-int-2}.
Using the same argument as the one introduced  in deriving
\eqnref{phi_blow_up_1}, one can show that the following lemma
holds.
\begin{lem}\label{represen}
Let $g$ and $\vp$ be the functions given in \eqnref{cd}. Then $g$
and $\vp$ are given by
 \begin{equation}\label{sr}
 \begin{array}{rcl}
 g &=& 2\ds\sum_{m=0}^{+ \infty}\ds\frac{1}{(2\lambda)^m}
  \ds\frac{\p}{\p\nu_\Om}[(R_B R_\Om)^m(I-\ds\frac{1}{2\lambda}R_B)
  \Dcal_\Om f] \quad \mbox{on } \p \Om, \\
 \nm \vp &=& \ds\frac{1}{\lambda}\sum_{m=0}^{+\infty}
  \ds\frac{1}{(2\lambda)^m}\ds\frac{\p}{\p\nu_B}
  [(R_\Om R_B)^m(I-R_\Om)\Dcal_\Om f] \quad \mbox{on } \p B,
 \end{array}
 \end{equation}
where the series in \eqnref{sr} converge absolutely and uniformly.
\end{lem}

\medskip
\noindent{\sl Proof of Lemma \ref{represen}}. The convergence of
the formula \eqnref{sr} will be proved in the course of proving
Theorem \ref{secondmain} in Section 4.

We first prove that for $(h_1, h_2) \in L^2_0(\p\Om) \times
L^2_0(\p B)$ there exists a unique solution $(g,\vp)\in
L^2_0(\p\Om) \times L^2_0(\p B)$ such that
$$
\left\{
\begin{array}{rcl}
\ds\frac{1}{2}g - \frac{\p (\Scal_B \vp)}{\p\nu_\Om} & = &\ds
h_1
\quad \mbox{on } \p \Om, \\
\nm
\ds \lambda\vp + \frac{\p (\Scal_\Om g)}{\p\nu_B}& =& \ds h_2 \quad \mbox{on
} \p B .
\end{array} \right.
$$
Since $B$ is away from $\p \Om$, the operator $(g, \vp) \to (-
\frac{\p (\Scal_B \vp)}{\p\nu_\Om}, \frac{\p (\Scal_\Om
g)}{\p\nu_B})$ is a compact operator on $L^2_0(\p\Om) \times
L^2_0(\p B)$. Thus, by the Fredholm alternative, it suffices to
prove that if $(h_1, h_2)=(0,0)$, then the solution $(g,
\vp)=(0,0)$. In order to do this, suppose that
$$
\left\{
\begin{array}{rcl}
\ds\frac{1}{2}g - \frac{\p (\Scal_B \vp)}{\p\nu_\Om} & = &0
\quad \mbox{on } \p \Om, \\
\nm
\ds \lambda\vp + \frac{\p (\Scal_\Om g)}{\p\nu_B}& =& 0 \quad \mbox{on } \p
B .
\end{array} \right.
$$ Then the function $u= \Scal_\Om (g) + \Scal_B(\vp)$
in $\Om$ is a solution of $\nabla \cdot (1+ (\frac{1}{k}-1)
\chi(B)) \nabla u=0$ in $\Om$ and satisfies $\frac{\partial
u}{\partial \nu} |_- =0$ on $\p\Om$. This implies that $u$ is
constant in $\Omega$, and hence
 $\Scal_B(\vp)$ is harmonic in $\Om$. It then follows
from the jump formula \eqnref{singlejump} that $\vp=0$ and
therefore, $ \Scal_\Om (g) = \mbox{constant}$ in $\Om$. Hence, $g=
- 2 \frac{\partial \Scal_\Om g}{\partial \nu}|_- =0$.

We now prove that the pair $(g, \vp)$ given by \eqnref{sr} satisfies
\eqnref{sy}. Observe that the function
 $$
 (R_B R_\Om)^m(I-\frac{1}{2\lambda}R_B) (\Dcal_\Om f)(X)= (\Dcal_\Om f)
((R_\Om R_B)^m (X))- \frac{1}{2\lambda} (\Dcal_\Om f) (R_B(R_\Om
R_B)^m (X))
 $$
is harmonic in $\RR^2 \setminus \overline B$ and approaches to
 $$
 (\Dcal_\Om f)
((R_\Om R_B)^{m-1} R_\Om (Z))- \frac{1}{2\lambda} (\Dcal_\Om f)
(R_B(R_\Om R_B)^{m-1} R_\Om (Z))
 $$
as $|X| \to +\infty$,  where $Z$ is the center of $B$. Since
$\Scal_\Om (1)$ is constant in $\Om$, it follows from
\eqnref{lm12} that
\begin{align*}
\pd{(\Scal_\Om g)}{\nu_B} & = \sum_{m=0}^{+\infty}
\frac{1}{(2\lambda)^m}
  \frac{\p}{\p\nu_B}[R_\Om (R_B R_\Om)^m(I-\ds\frac{1}{2\lambda}R_B)
  \Dcal_\Om f] \\
&= \sum_{m=0}^{+\infty} \frac{1}{(2\lambda)^m}
  \frac{\p}{\p\nu_B} [(R_\Om R_B)^m (R_\Om - I)
  \Dcal_\Om f] + \pd{(\Dcal_\Om f)}{\nu_B}.
\end{align*}
Likewise one can show that
$$
\pd{(\Scal_B \vp)}{\nu_\Om} =
\ds\sum_{m=0}^{+\infty}\ds\frac{1}{(2\lambda)^m}
  \ds\frac{\p}{\p\nu_\Om}[(R_B R_\Om)^m(I-\ds\frac{1}{2\lambda}R_B)
  \Dcal_\Om f] - \pd{(\Dcal_\Om f)}{\nu_\Om}.
$$
Thus $(g, \vp)$ satisfies \eqnref{sy} and the proof is complete.
\qed

The representation of $g$ and $\vp$ given in \eqnref{sr} can be
simplified using the relation
 \be \label{reflec}
 R_\Om\Dcal_\Om f (X) =\Dcal_\Om f (R_\Om (X))= - \Dcal_\Om f(X) +
 \mbox{constant}, \quad X \in \RR^2 \setminus \p\Om,
 \ee
which follows from \eqnref{doublejump} and \eqnref{KDKD} since
$\Om$ is a disk. Using \eqnref{reflec}, we then compute
 \begin{align}
 g &= 2\ds\sum_{m=0}^{+\infty}\ds\frac{1}{(2\lambda)^m}
  \ds\frac{\p}{\p\nu_\Om}[(R_B R_\Om)^m(I-\ds\frac{1}{2\lambda}R_B)
  \Dcal_\Om f]  \nonumber \\
 & = 2 \sum_{m=0}^{+\infty}\ds\frac{1}{(2\lambda)^m}
  \ds\frac{\p}{\p\nu_\Om}[(R_B R_\Om)^m(I+ \frac{1}{2\lambda}R_B R_\Om)
  \Dcal_\Om f] \nonumber \\
 & = 4 \sum_{m=1}^{+\infty} \frac{1}{(2\lambda)^m}
  \ds\frac{\p}{\p\nu_\Om}[(R_B R_\Om)^m
  \Dcal_\Om f] + 2 \frac{\p}{\p\nu_\Om} \Dcal_\Om f \quad \mbox{on } \p\Om. \label{srgnew}
 \end{align}
Likewise we can show that
 \be \label{srvpnew}
  \vp = \frac{2}{\lambda}
\sum_{m=0}^{+\infty} \frac{1}{(2\lambda)^m}\ds\frac{\p}{\p\nu_B}
  [(R_\Om R_B)^m\Dcal_\Om f] \quad \mbox{on } \p B.
 \ee

The following lemma is also of importance to us.
\begin{lem}\label{conj-2}
Let $u$ be the solution of \eqnref{eq-neu} for $g \in
\mathcal{C}^\alpha(\p \Om)$ and let $G$ be the function satisfying
$\pd{G}{T}=g$ on $\p\Om$ and $\int_{\p\Om} G =0$. Define $v$ to be
the solution of the following conductivity problem
\begin{equation} \label{eqv1}
\ \left \{
\begin{array}{l}
\ds \nabla  \cdot \Bigr( 1 + (\frac{1}{k} -1) \chi(B) \Bigr)
\nabla v =0 \quad \mbox{in } \Om, \\ \nm \ds v=G \quad\mbox{on }~
\p\Om.
\end{array}
\right .
\end{equation}
Then
\be \label{onB}
\frac{\partial u}{\partial  T} =  -\frac{\partial v}{
\partial \nu} \bigg |_+  \quad \mbox{ on }~ \p B.
\ee
Moreover, $\Dcal_\Om (v|_{\p\Om})$ is a harmonic conjugate to
$\Scal_\Om g$ in $\Om$.
\end{lem}

\pf Let $w$ be a harmonic conjugate of $u$ in $\Om \setminus
\overline B$ and $B$. Such a conjugate function exists in $\Om
\setminus \overline B$ since $\int_{C} \pd{u}{\nu} d\sigma=0$ for
any simple closed curve $C$ in $\Om \setminus \overline B$.
Moreover, since $u$ is $\mathcal{C}^{1,\alpha}$, so is $w$.

Define $v$ by \be \label{defv} v(X) := \begin{cases}
\ds w(X), \quad & X \in \Om \setminus \overline{B}, \\
\ds kw(X) - \frac{k}{|\partial B|} \int_{\partial B} w\; d\sigma,
& X \in B.
\end{cases}
\ee Then one can see from the Cauchy-Riemann equation and the
transmission conditions on $u$ that
\begin{align*}
\pd{v}{T} \bigg|_{+} =  \pd{u}{\nu} \bigg |_{+} = k \pd{u}{\nu} \bigg
|_{-} = \pd{v}{T} \bigg|_{-}, \\
\pd{v}{\nu} \bigg|_{+} = -\pd{u}{T} \bigg|_{+} = - \pd{u}{T}
\bigg|_{-} = \frac{1}{k} \pd{v}{\nu} \bigg|_{-}.
\end{align*}
Thus $v$ defined by \eqnref{defv} is the unique solution to
\eqnref{eqv}, and hence \eqnref{onB} holds.

It follows from \eqnref{doublejump} that
$$
\Dcal_\Om (v|_{\p\Om})|_{-} = \frac{1}{2} v \quad \mbox{on } \p\Om,
$$
and hence
$$
\pd{(\Dcal_\Om (v|_{\p\Om}))}{T} = \frac{1}{2} \pd{v}{T} =\frac{1}{2} g = -
\pd{(\Scal_\Om g)}{\nu}
  \quad \mbox{on } \p\Om.
$$
Therefore $\Dcal_\Om (v|_{\p\Om})$ is a harmonic conjugate of
$\Scal_\Om g$ in $\Om$. This completes the proof. \qed

\section{Proof of Theorem \ref{firstmain}}\label{sec1}
 At this point we have all the necessary ingredients to prove
Theorem \ref{firstmain}. As has been said, the lower bound in
\eqnref{estthm} and \eqnref{estthm1} is an improvement of the one
obtained in \cite{AKL}.

Recall that there are two disks $B_j=B(Z_j, r_j)$, $j=1,2$, inside
$\Om$, and that $R_j$ is the reflection with respect to $\p B_j$.
We suppose that both centers $Z_1$ and $Z_2$ are on the $x$-axis.

If $X$ is on the $x-$axis, {\em i.e.}, $X=(x,0)$, straightforward
calculations show that $$ D R_{i} (X)=g_i (X) \left (
\begin{array}{cc}
  -1& 0 \\
  0 & 1
\end{array}
\right ),
$$
and
\be \label{rif}
\nabla (R_i f)(X) =  \nabla f (R_i(X)) g_i (X)\left (
\begin{array}{cc}
-1& 0 \\
  0 & 1
\end{array}
\right ), \ee where \be \label{gix}
g_i(X)=\frac{r_i^2}{|X-Z_i|^2}, \quad i=1,2. \ee Therefore,
\begin{equation} \label{phi_blow_up_2}
\nabla ((R_2R_1)^m H)(X ) = \bigg[\prod_{i=1}^{2m}
g_{l_i}(R_{l_{i-1}} \cdots R_{l_1} (X ))\bigg]
{\nabla H((R_1 R_2)^m (X ))} \begin{pmatrix}
  1 & 0 \\
  0 & 1
\end{pmatrix},
\end{equation}
and
\begin{equation}\label{phi_blow_up_3}
\begin{array}{l}
\ds \nabla ((R_2R_1)^m R_2 H)(X ) \\
\nm \ds  =  g_2((R_1R_2)^m(X ))
\bigg[\prod_{i=1}^{2m} g_{l_i}(R_{l_{i-1}} \cdots R_{l_1} (X
))\bigg] \nabla H(R_2(R_1 R_2)^m (X ))
\begin{pmatrix}
-1 & 0 \\
  0 & 1
\end{pmatrix},
\end{array}
\end{equation}
where $ (l_1,\dots,l_{2m})=
(\overbrace{2,1,2,1,\dots,2,1}^{m\mbox{ times}})$.

Let $u$ be the solution to the equation \eqnref{cond-trans}.
Combining \eqnref{singlejump}, \eqnref{KDKD}, \eqnref{repp}, and
\eqnref{multi-int-2} yields
$$
\pd{u}{\nu^{(i)}} \bigg |_{\pm}
=\pd{H}{\nu^{(i)}}+\pd{(\Scal_{B_2}\varphi_2)}{\nu^{(i)}}
\bigg |_{\pm}   + \pd{(\Scal_{B_1}\varphi_1)}{\nu^{(i)}} \bigg |_{\pm}
=(\lambda_i \pm \frac{1}{2})\varphi_i \quad \mbox{on } \partial
B_i, \ i=1,2.
$$
Consequently,
\begin{equation} \label{eqpm}
| \nabla u |_{\pm}(X_i) | \ge \left |
\pd{u}{\nu^{(i)}}\big|_\pm(X_i) \right | \geq \left | \lambda_i
\pm \frac{1}{2} \right | | \varphi_i(X_i) |. \ee

Suppose that $k_1>1$ and $k_2>1$. By \eqnref{phi_blow_up_1},
\eqnref{phi_blow_up_2} and \eqnref{phi_blow_up_3}, we obtain the
following inequality:
\begin{align*}
|\varphi_1(X_1)|&\geq  \frac{1}{\lambda_1}\ds\sum_{m=0}^{+\infty}
\frac{1}{(4\lambda_1\lambda_2)^m}
\bigg(a^{2m}+\frac{1}{2\lambda_2}a^{2m+1}  \bigg)\ds\inf_{I}
\big|\nabla H\cdot\nu^{(1)} \big|\\& \geq
\frac{1}{\lambda_1}\frac{1+\frac{a}{2\lambda_2}}{1-\frac{a^2}{4\lambda_1\lambda_2}
}\ds\inf_{I} \big|\nabla H\cdot\nu^{(1)} \big|, \end{align*} where
$a := ({1 + 2 (r_*/r_{\min})\sqrt{\epsilon}})^{-1}$.  We then get
from \eqnref{eqpm} that
\begin{align*}
|\nabla u |_+ (X_1)| \geq\frac{ C}{1 -\tau + (r_*/r_{\min})
\sqrt{\epsilon}}\ds\inf_{I} \big|\nabla H\cdot\nu^{(1)} \big|.
\end{align*}
Likewise we obtain
\begin{align*}
|\nabla u |_+ (X_2)| \geq\frac{ C}{1 -\tau + (r_*/r_{\min})
\sqrt{\epsilon}}\ds\inf_{I} \big|\nabla H\cdot\nu^{(2)} \big|.
\end{align*}
Thus the lower bound in \eqnref{estthm} is now derived.

If both $k_1$ and $k_2$ are less than 1, then one can use Lemma
\ref{conj} to obtain the lower bound in \eqnref{estthm1}. See
\cite{AKL} for the details.

\bigskip
We now prove \eqnref{estthm-upper}. We need the following lemmas.

\begin{lem} \label{lem31}
Let $r_{\max}, r_{\min}, r_*,$ and $\epsilon$ be as in Theorem
\ref{firstmain}. If $\epsilon$ is small enough, then for any $X
\in \overline{B}_1$ and $n \geq 8 r_*/ \sqrt{\epsilon}$ we have
$$
|(R_1R_2)^n(X)- Z_2|\geq r_2 (1 + \frac{r_*}{2r_{\max}}
\sqrt{\epsilon}),
$$
and
$$
|(R_2R_1)^n R_2(X)- Z_1|\geq r_1 (1 + \frac{r_*}{2r_{\max}}
\sqrt{\epsilon}).
$$
For any $X \in \overline{B}_2$ and  $n \geq 8 r_*/ \sqrt{\epsilon}$
we have
$$
|(R_2R_1)^n(X)- Z_1|\geq r_1 (1 + \frac{r_*}{2r_{\max}}
\sqrt{\epsilon}),
$$
and
$$
|(R_1R_2)^n R_1(X)- Z_2|\geq r_2 (1 + \frac{r_*}{2r_{\max}}
\sqrt{\epsilon}).
$$
\end{lem}

\pf After a translation and a rotation if necessary, we may assume
that $B_1 = B((0,0), r_1)$ and $B_2 = B((r_1+r_2+\epsilon,
0),r_2)$, {\em i.e.}, $Z_1=(0,0)$ and $Z_2= (r_1+r_2+\epsilon,
0)$. It is easy to show that the fixed points of the combined
reflections $R_2R_1$ and $R_1R_2$ are the points $(x_i, 0)$, for
$i=1,2$, where $x_1$ and $x_2$ are the roots of the quadratic
equation
$$
\bigg(r_1+r_2+\epsilon \bigg) x^2 + \bigg(r_2^2 -r_1^2 -(r_1
+r_2+\epsilon)^2 \bigg) x + r_1^2 \bigg(r_1+r_2+\epsilon \bigg) = 0.
$$
Then, as $\epsilon$ goes to zero,
$$
x_1= r_1{-} \sqrt{\frac{2r_1r_2}{r_1+r_2}} \sqrt{\epsilon} + O(\epsilon),
\quad x_2 = r_1{+} \sqrt{\frac{2r_1r_2}{r_1+r_2}} \sqrt{\epsilon} +
O(\epsilon),
$$
and \be \label{rstar}\frac{1}{2} r_* \sqrt{\epsilon}\leq | r_1-x_j
| \leq 2r_* \sqrt{\epsilon}, j=1,2. \ee

Let $X_1=(r_1,0)$ and $(t_n,0)=(R_1R_2)^n(X_1)$. We have
$$
t_{n+1}=\frac{r_1^2}{(r_1+r_2+\epsilon)-
\ds \frac{r_2^2}{r_1+r_2+\epsilon-t_n}},
$$
and hence
\begin{align*}
| t_{n+1}-x_1 | &= \left | \frac{r_1^2}{(r_1+r_2+\epsilon)-
\ds \frac{r_2^2}{r_1+r_2+\epsilon-t_n}}
-\frac{r_1^2}{(r_1+r_2+\epsilon)-
\ds \frac{r_2^2}{r_1+r_2+\epsilon-x_1}} \right | \\
&=|t_n-x_1|
\left |
\frac{r_1r_2}{(r_1+r_2+\epsilon)(r_1+r_2+\epsilon-t_n)-r_2^2}\right|
\\
& \quad\quad  \times \left |
\frac{r_1r_2}{(r_1+r_2+\epsilon)(r_1+r_2+\epsilon-x_1)-r_2^2}\right| \\
&\leq \frac{|t_n-x_1|}{1+\frac{\sqrt{\epsilon}}{r_*}}.
\end{align*}
Thus, \be \label{tnx} |t_n-x_1|\leq \frac{r_*}{4} \sqrt{\epsilon}
\quad\mbox{if }     n\geq \frac{8r_*}{\sqrt{\epsilon}}. \ee
Observe that $|(R_1R_2)^n(X_1)-Z_2| = r_1 - t_n + r_2 + \epsilon$.
Therefore, it follows from \eqnref{rstar} and \eqnref{tnx} that
\begin{equation}\label{fixedd}
|(R_1R_2)^n(X_1)-Z_2|\geq
r_2+\frac{r_*}{2}\sqrt{\epsilon} \quad \mbox{if } n\geq
\frac{8r_*}{\sqrt{\epsilon}}.
\end{equation}

Let $X\in \overline B_1$. For any $0\leq t\leq r_1$ such that
$|X-Z_2|\geq |(t,0)-Z_2|$, we can easily see that
$$
|R_2(X)-Z_1|\geq |R_2(t,0)-Z_1| \quad \mbox{and} \quad
|R_1R_2(X)-Z_2|\geq |R_1R_2(t,0)-Z_2|.
$$
Since $R_1R_2(t,0)=(s,0)$ for some $s$ satisfying $0\leq s\leq r_1$,
then, by repeating the
above inequalities, we have for each positive integer $m$,
\begin{align*}
|(R_1R_2)^m(X)-Z_2| &\geq |(R_1R_2)^m(t,0)-Z_2|, \\
\nm
|R_2(R_1R_2)^m(X)-Z_1| &\geq |R_2(R_1R_2)^m(t,0)-Z_1|.
\end{align*}
In particular, for the case $t=r_1$, we obtain that
\begin{align*}
|(R_1R_2)^m(X)-Z_2| & \geq |(R_1R_2)^m(X_1)-Z_2|, \\
\nm
|R_2(R_1R_2)^m(X)-Z_1| &\geq |R_2(R_1R_2)^m(X_1)-Z_1|, \quad \forall
\; X \in \overline{B}_1.
\end{align*}
Similarly, we have
\begin{align*}
|(R_2R_1)^m(X)-Z_1| & \geq |(R_2R_1)^m(X_2)-Z_1|, \\
\nm
|R_1(R_2R_1)^m(X)-Z_2| &\geq |R_1(R_2R_1)^m(X_2)-Z_2|, \quad \forall
\; X \in \overline{B}_2.
\end{align*}
Therefore the first two inequalities in Lemma \ref{lem31} follow
from \eqnref{fixedd}. The second pair of inequalities can be
derived by interchanging $B_1$ and $B_2$. This completes the
proof. \qed

\begin{lem}
\begin{itemize}
\item[{\rm (i)}]  Let $X_1=(r_1,0)$.
We have
\begin{equation} \label{gx1}
\ds g_2((R_1R_2)^m(X_1)), \ g_1((R_2R_1)^m R_2(X_1)) \geq
\frac{1}{1 + 8 (r_*/r_{\min}) \sqrt{\epsilon}},  \quad \forall\;
m\in \NN.
\end{equation}
\item[{\rm (ii)}]
For all $X\in \overline{B}_1$, we have
\begin{align}
\begin{cases} \label{gx2}
\ds g_2((R_1R_2)^m(X))\leq 1,  \quad \forall\; m\in \NN,\\
\nm \ds g_2((R_1R_2)^m(X))\leq \frac{1}{1 +
  (r_*/r_{\max}) \sqrt{\epsilon}}, \quad \quad \forall\;  m \geq
  8 r_{*}/\sqrt{\epsilon},
\end{cases}
\end{align}
and similarly, for all $X\in \overline{B}_2$, we have
\begin{align}
\begin{cases} \label{gx3}
\ds g_1((R_2R_1)^m(X))\leq 1, \quad \forall\; m\in \NN,\\
\nm \ds g_1((R_2R_1)^m(X))\leq \frac{1}{1 +
  (r_*/r_{\max})\sqrt{\epsilon}},  \quad \quad \forall\; m \geq
8 r_{*}/\sqrt{\epsilon}.
\end{cases}
\end{align}
\end{itemize}
\end{lem}

\pf Since $(R_1R_2)^m(X_1)$ is between $\widetilde{X}_1$ and
$\widetilde{X}_2$ where $\widetilde{X}_1 =(x_1,0)$ and
$\widetilde{X}_2=(x_2,0)$
$(x_1<x_2)$
are fixed points of $R_2R_1$ and $R_1R_2$, respectively, we have
\begin{align*}
g_2((R_1R_2)^m(X_1))&=\frac{r_2^2}{|(R_1R_2)^m(X_1)-Z_2|^2}\geq
\frac{r_2^2}{|\widetilde{X}_1-Z_2|^2}\\&=\frac{r_2^2}{(r_1+r_2+\epsilon-x_1)^2}\geq
\frac{r_2^2}{(r_2+2r_{*}\sqrt{\epsilon}+\epsilon)^2}\geq
\frac{1}{1 + 8 (r_*/r_{\min}) \sqrt{\epsilon}}.
\end{align*}
The second inequality in \eqnref{gx1} can be proved in exactly the
same way. Lemma \ref{lem31} and the definition of $g_i$ give the
upper bounds \eqnref{gx2} and \eqnref{gx3}. \qed

To establish our upper bound we first observe that since
$u(X)-H(X) \to 0$ as $|X| \to + \infty$, $|\nabla(u-H)|$ attains
its maximum on either $\p B_1$ or $\p B_2$, and hence
\begin{align}
\| \nabla u \|_{L^\infty(\Om\setminus \overline{B_1 \cup B_2})} &
\le \| \nabla (u-H) \|_{L^\infty(\Om\setminus \overline{B_1 \cup
B_2})}
+ \| \nabla H \|_{L^\infty(\Om)} \nonumber \\
& \le \| \nabla (u-H)|_{+} \|_{L^\infty(\p B_1 \cup \p B_2)}
+ \| \nabla H \|_{L^\infty(\Om)} \nonumber \\
& \le \| \nabla u|_{+} \|_{L^\infty(\p B_1 \cup \p B_2)} + \|
\nabla H \|_{L^\infty(\Om)} , \label{linfty}
\end{align}
and
 $$
 \| \nabla u \|_{L^\infty(B_1 \cup B_2)} \le \| \nabla u|_{-} \|_{L^\infty(\p B_1 \cup \p B_2)}.
 $$
We also have from \eqnref{repp} and \eqnref{eqpm} that
 \be \label{uLi}
 \|\nabla u|_{\pm}\|_{L^\infty(\p B_i)}\leq  (|\lambda_i| + \frac{1}{2} )
 \;
  \|  \varphi_i   \|_{L^\infty(\p B_i)} +\bigg\|\pd{u}{T}
 \bigg\|_{L^\infty(\p B_i)}.
 \ee

Let $N$ be the first integer larger than $8 r_*/\sqrt{\ep}$. It
then follows from \eqnref{phi_blow_up_1}, \eqnref{gx2}, and
\eqnref{gx3} that
 \begin{align*}
 | \varphi_1(X ) | \ds & \leq  \;
 \|\nabla H \|_{L^\infty(B_1 \cup B_2)} \frac{1}{|\lambda_1|} \bigg( \ds
 \sum_{m < N}
 \frac{1}{|4\lambda_1\lambda_2|^m}
 \big(1+\frac{1}{2|\lambda_2|}\big ) \\
& \quad + \frac{1}{|4\lambda_1\lambda_2|^N} \sum_{m=0}^{+\infty}
\frac{1}{|4\lambda_1\lambda_2|^m}
\big(b^{2m}+\frac{1}{2|\lambda_2|}b^{2m+1}\big)  \bigg),
\end{align*}
for any $X\in \p B_1$, where $b := {1}/({1 +
  (r_*/r_{\max})\sqrt{\epsilon}})$. Thus, for each $X \in \p B_1$,
\begin{align*}
| \varphi_1(X ) | &\leq \frac{C\|\nabla H \|_{L^\infty(B_1 \cup
B_2)}}{|\lambda_1|} \bigg(\frac{1}{1 -|\tau| + r_*/r_{\max}
\sqrt{\epsilon}}  + \frac{1-
|\tau|^{8r_*/\sqrt{\epsilon}}}{1-|\tau|}\bigg)\\&\leq
C\frac{\|\nabla H \|_{L^\infty(B_1 \cup B_2)}}{|\lambda_1|(1
-|\tau| + r_*/r_{\max} \sqrt{\epsilon})},
\end{align*} for any for any $X\in \p B_1$.
Similarly, we have the desired estimate for $\vp_2(X)$, $X\in \p
B_2$ and hence, \be \label{pdpu}
\bigg\|\pd{u}{\nu}\bigg\|_{L^\infty(\p B_1 \cup \p B_2)}\leq
\frac{C\|\nabla H \|_{L^\infty(B_1 \cup B_2)}}{1 -|\tau| +
r_*/r_{\max} \sqrt{\epsilon}}. \ee

To estimate $\p u/\p T$ we use Lemma \ref{conj}. Let $\widetilde
H$ be a harmonic conjugate of $H$ and $v$ be the solution to
\eqnref{eqv}. Then by \eqnref{pdpu}, we have
$$
\bigg\|\pd{v}{\nu}\bigg\|_{L^\infty(\p B_1 \cup \p B_2)}\leq
\frac{C\|\nabla \widetilde H \|_{L^\infty(B_1 \cup B_2)}}{1
-|\tau| + r_*/r_{\max} \sqrt{\epsilon}}.
$$
Since $\|\nabla \widetilde H \|_{L^\infty(B_1 \cup B_2)} \|\nabla
H \|_{L^\infty(B_1 \cup B_2)}$, it follows from Lemma \ref{conj}
that \be \label{pdpuT} \bigg\|\pd{u}{T}\bigg\|_{L^\infty(\p B_1
\cup \p B_2)}\leq \frac{C\|\nabla H \|_{L^\infty(B_1 \cup B_2)}}{1
-|\tau| + r_*/r_{\max} \sqrt{\epsilon}}. \ee Combining
\eqnref{linfty}, \eqnref{uLi}, \eqnref{pdpu}, and \eqnref{pdpuT}
yields the upper bound in \eqnref{estthm-upper} and completes the
proof.

\section{Proofs of Theorems \ref{secondmain} and \ref{thirdmain}}

We suppose that $\Omega=B((0,0),\rho)$ and
$B=B((\rho-r-\epsilon,0),r)$ after rotation and translation if
necessary, and that $\epsilon\ll \rho-r$. The conductivities of
$\Omega$ and $B$ are $1$ and $k$, $0<k\neq1<+\infty$,
respectively. Let
$$
g_\Om(X)=\ds\frac{\rho^2}{|X|^2}, \quad
g_B(X)=\ds\frac{r^2}{|X-(\rho-r-\epsilon,0)|^2}.
$$
The functions $g_\Om$ and $g_B$ play the roles of $g_1$ and $g_2$
in the previous section. Let $R_\Om$ and $R_B$ be the reflections
with respect to $\Om$ and $B$.

\begin{lem}\label{lem303}
Let $X_1=(\rho-\epsilon,0)$, the point on $\overline{B}$ closest
to $\p\Om$. For any positive integer $n$ and $X\in\overline{B}$,
we have
 \begin{equation} \label{gBRom}
 g_B(R_\Omega(R_B R_\Omega)^n(X))g_\Omega((R_B R_\Omega)^n(X))\leq
 g_B(R_\Omega(R_B R_\Omega)^n(X_1))g_\Omega((R_B R_\Omega)^n(X_1)).
 \end{equation}
\end{lem}

\pf For any $X=(x,y)$, one can easily see that
\begin{align}\label{ew}
g_B(R_\Omega(X))g_\Omega(X)&=\frac{r^2}{\big(
\frac{\rho^2x}{x^2+y^2}-(\rho-r-\epsilon)\big)^2+\big(\frac{\rho^2y}{x^2+y^2}\big)^2}\cdot\frac{\rho^2}{x^2+y^2}\nonumber\\
&=\frac{\rho^2
r^2}{(\rho-r-\epsilon)^2}\cdot\frac{1}{\big(\frac{\rho^2}{\rho-r-\epsilon}-x\big)^2+y^2}.
\end{align}
Since $\frac{\rho^2}{\rho-r-\epsilon} > \rho -\epsilon$, it
immediately follows that
$$
g_B(R_\Omega(X))g_\Omega(X)\leq
g_B(R_\Omega(X_1))g_\Omega(X_1), \quad \forall X\in\overline{B}.
$$
If $X$ satisfies
$|X-(\rho-r-\epsilon)|\leq|(t,0)-(\rho-r-\epsilon)|$, with $t >
\rho-r-\epsilon$,  then
$$|R_B R_\Omega(X)-(\rho-r-\epsilon)|\leq|R_B
R_\Omega(t,0)-(\rho-r-\epsilon)|.$$ Using this fact repeatedly, we
have
$$
|(R_B R_\Omega)^n(X)-(\rho-r-\epsilon)|\leq|(R_B
R_\Omega)^n(X_1)-(\rho-r-\epsilon)|.
$$
By combining this inequality with \eqnref{ew}, we obtain that, for
$X \in \overline{B}$,
\begin{equation*} g_B(R_\Omega(R_B R_\Omega)^n(X))g_\Omega((R_B
R_\Omega)^n(X))\leq g_B(R_\Omega(R_B
R_\Omega)^n(X_1))g_\Omega((R_B R_\Omega)^n(X_1)),
\end{equation*}
which completes the proof. \qed

Recall that $P_1$ and $P_2$ are the fixed points of the combined
reflections $R_BR_\Om$ and $R_\Om R_B$, respectively. Observe that
$P_1 \in B$. If $P_i=(x_i,0)$ for $i=1,2$, then $x_i \ (x_1<x_2)$
are the roots of the quadratic equation
 $$
 (\rho-r-\epsilon)x^2+
 (r^2-\rho^2-(\rho-r-\epsilon)^2)x+\rho^2(\rho-r-\epsilon)=0.
 $$
It then follows that
 \be \label{x1x2}
 x_1= \rho{-} \sqrt{\frac{2\rho r}{\rho-r}} \sqrt{\epsilon} +
 O(\epsilon), \quad x_2 = \rho{+} \sqrt{\frac{2\rho r}{\rho-r}}
 \sqrt{\epsilon} + O(\epsilon).
 \ee
Moreover,
$$ \frac{\sqrt{\epsilon}}{r_*} \leq | x_j - \rho| \leq 2
\frac{\sqrt{\epsilon}}{r_*}, j=1,2, $$ for $\epsilon$ small enough.

As a direct consequence of \eqnref{gBRom} and \eqnref{x1x2}, we
have the following lemma which plays a crucial role in deriving
the lower bound \eqnref{est}.

\begin{lem} \label{lem402}
For each positive integer $n$, the following inequality holds:
 \be \label{gBR}
g_B(R_\Om(R_BR_\Om)^n(X_1)) g_\Om((R_BR_\Om)^n(X_1)) \geq
\frac{1}{1+4 r^*\sqrt{\epsilon}}. \ee
\end{lem}

\pf Since $P_1 \in B$, we have
 \begin{align*}
 g_B(R_\Om(R_BR_\Om)^n(X_1))
 g_\Om((R_BR_\Om)^n(X_1)) & \geq g_B(R_\Om(R_BR_\Om)^n(P_1))
 g_\Om((R_BR_\Om)^n(P_1)) \\
 & \geq g_B(R_\Om(P_1))g_\Om(P_1).
 \end{align*}
On the other hand, it follows from \eqnref{ew} and \eqnref{x1x2}
that
 $$
 g_B(R_\Omega(P_1))g_\Omega(P_1) =\frac{\rho^2
 r^2}{\big(\rho^2-x_1 (\rho-r-\epsilon) \big)^2}  \geq
 \frac{1}{1+4 r^*\sqrt{\epsilon}}.
 $$
The proof is complete. \qed

The following lemma is also of use to us.
\begin{lem} \label{lem403}
For each positive integer $n$, similarly to \eqnref{gBR} the
following holds:
 \be\label{estgbgom}
 g_\Om(R_B (R_\Omega R_B)^{n}(X_2)) g_B((R_\Om R_B)^{n}(X_2))
 \geq
 \frac{1}{1+4r^*\sqrt{\epsilon}}.
 \ee
\end{lem}

\pf Since the proof of Lemma \ref{lem403} is parallel to that of
Lemma \ref{lem402}, we very briefly sketch it. As before, since
$P_2 \in \RR^2 \setminus \Om$, we can show that
 \begin{align*}
  g_B((R_\Om R_B)^{n}(X_2))
 g_\Om((R_BR_\Om)^{n}R_B(X_2)) & \ge  g_B((R_\Om R_B)^{n}(P_2))
 g_\Om((R_BR_\Om)^{n}R_B(P_2)) \\
 & = g_B(P_2) g_\Om(R_B(P_2)).
 \end{align*}
Since, as one can see easily, $R_\Om (P_1) = P_2$ we have
 $$
 g_B(P_2) g_\Om(R_B(P_2)) = g_B(R_\Om(P_1)) g_\Om(P_1),
 $$
and hence \eqnref{estgbgom} follows. \qed

 Next, we need the following lemma to derive
the upper bound.
\begin{lem}\label{lem304}
For each positive integer $n$,
\begin{equation} \label{lessone}
g_B(R_\Omega(R_B R_\Omega)^n(X))g_\Omega((R_B R_\Omega)^n(X))\leq
1, \quad X \in \overline{B},
\end{equation}
and for $n\geq \frac{1}{4 r^* \sqrt{\epsilon}},$
\begin{equation} \label{lesstwo}
g_B(R_\Omega(R_B R_\Omega)^n(X))g_\Omega((R_B R_\Omega)^n(X))\leq
\frac{1}{1+r^*\sqrt{\epsilon}}, \quad X \in \overline{B}.
\end{equation}
\end{lem}

\pf The inequality \eqnref{lessone} is obvious. We only have to
prove \eqnref{lesstwo}. Let $(t_n,0):=(R_BR_\Omega)^n(X_1)$. Then
we have
$$
t_{n+1}=\frac{r^2}{
\frac{\rho^2}{t_n}-(\rho-r-\epsilon)}+ \rho-r-\epsilon.
$$
Recall that $P_1=(x_1, 0)$. We then have
\begin{align*}
| t_{n+1}-x_1 | &= \bigg| \frac{r^2}{
\frac{\rho^2}{t_n}-(\rho-r-\epsilon)}-\frac{r^2}{
\frac{\rho^2}{x_1}-(\rho-r-\epsilon)} \bigg|\\
&=|t_n-x_1| \left[\frac{\rho
r}{\rho^2-(\rho-r-\epsilon)x_1}\right]
\left[\frac{\rho r}{\rho^2-(\rho-r-\epsilon)t_n}\right]\\
&\leq \frac{|t_n-x_1|}{1+\sqrt{\frac{\rho-r}{\rho
r}}\sqrt{\epsilon}}.
\end{align*}
If $n\geq \frac{1}{4 r^*\sqrt{\epsilon}}$, we get
 $$
 |t_n-x_1|\leq \frac{1}{2}\sqrt{\frac{\rho r}{\rho-r}}\sqrt{\epsilon},
 $$
and therefore
\begin{equation}\label{fixed}
|(R_B R_\Om)^n(X_1)|\leq
\rho{-} \sqrt{\frac{\rho r}{2(\rho-r)}} \sqrt{\epsilon}.
\end{equation}
By Lemma \ref{lem303} and \eqnref{ew},
\begin{align*}
g_B(R_\Omega(R_B R_\Omega)^n(X))g_\Omega((R_B R_\Omega)^n(X))&\leq
\frac{\rho^2r^2}{\Big(\rho^2-(\rho-r-\epsilon)\big(\rho-\sqrt{\frac{\rho
r}{2(\rho-r)}}\sqrt{\epsilon}\big)\Big)^2}\\
&\leq\frac{1}{1+r^*\sqrt{\epsilon}},
\end{align*}
which completes the proof. \qed

Since $R_B R_\Om(X)$ is in $\overline B$ for any
$X\in\overline\Om$, we immediately obtain the following corollary.

\begin{cor} \label{cor45}
For each positive integer $n$,
\begin{equation*}
g_B(R_\Omega(R_B R_\Omega)^n(X))g_\Omega((R_B R_\Omega)^n(X))\leq
1, \quad X \in \overline{\Om},
\end{equation*}
and for $n > \frac{1}{4 r^* \sqrt{\epsilon}},$
\begin{equation*}
g_B(R_\Omega(R_B R_\Omega)^n(X))g_\Omega((R_B R_\Omega)^n(X))\leq
\frac{1}{1+r^*\sqrt{\epsilon}}, \quad X \in \overline{\Om}.
\end{equation*}
\end{cor}

We are now ready to prove Theorem \ref{secondmain} and Theorem
\ref{thirdmain}.

\medskip
\noindent{\sl Proof of Theorem \ref{secondmain}}. Straightforward
computations yield that for $X_1=(\rho-\ep,0)$
\begin{align}
&\nabla ((R_\Om R_B)^m \Dcal_\Om f)(X_1) \nonumber \\
&  =\prod^{m-1}_{n=0}g_B(R_\Om(R_BR_\Om)^n(X_1))
g_\Om((R_BR_\Om)^n(X_1)) \nabla \Dcal_\Om f((R_B R_\Om)^m(X_1)) .
\label{nablaR}
\end{align}
Since $R_\Om(P_1)=P_2$ and $R_B(P_2)=P_1$, $(R_BR_\Om)^n(X_1)$
lies in $J_1$, the line segment between $P_1$ and $X_1$, for each
$n$. We may assume
 $$
 \inf_{X \in J_1} |\la \nabla
\Dcal_\Om f(X), \nu_B(X_1) \ra| \neq 0,
 $$
since otherwise the estimate \eqnref{est} is trivial. If $\ep$ is
small enough, then the length of $J_1$ is small and hence we may
further suppose that
 $$
 \la \nabla
 \Dcal_\Om f(X), \nu_B(X_1)\ra = \frac{\p (\Dcal_\Om f)}{\p
x}(X)
 $$
has the same sign for all $X \in J_1$. It then follows from
\eqnref{gBR} and \eqnref{nablaR} that
\begin{align}
&\bigg| \frac{\p}{\p\nu_B} [(R_\Om R_B)^m\Dcal_\Om f](X_1) \bigg| \nonumber \\
&  =\bigg| \prod^{m-1}_{n=0}g_B(R_\Om(R_BR_\Om)^n(X_1))
g_\Om((R_BR_\Om)^n(X_1)) \la \nabla \Dcal_\Om f((R_B
R_\Om)^m(X_1)), \nu_B(X_1) \ra \bigg| \nonumber \\
& \ge a^m \inf_{X \in J_1}\Big|\frac{\p (\Dcal_\Om f)}{\p
x}(X)\Big|, \label{romrb}
\end{align}
where $a:=(1+4 r^*\sqrt{\epsilon})^{-1}$.

Suppose that $k>1$ and hence $\lambda >0$. It follows from
\eqnref{romrb} that
\begin{align*}
  \Big| \vp(X_1)\Big|  &= \Big|\ds\frac{2}{\lambda}\sum_{m=0}^{+\infty}
  \ds\frac{1}{(2\lambda)^m}\ds\frac{\p}{\p\nu_B}
  [(R_\Om R_B)^m\Dcal_\Om f](X_1)\Big|\\
  &\geq\Big|\frac{2}{\lambda} \sum_{m=0}^{+\infty}
     \Big(\frac{a}{2\lambda}\Big)^m\Big|
     \cdot \inf_{X\in J_1}\Big|\frac{\p (\Dcal_\Om f)}{\p x}(X)\Big| \\
   &\geq\ds\frac{C \sigma \inf_{X\in J_1}\Big|\frac{\p (\Dcal_\Om f)}{\p
x}(X)\Big|}{1-\sigma + 4 r^*\sqrt{\ep}},
\end{align*}
since $\sigma=\frac{1}{2\lambda}$. Here $C$ is a constant which is
independent of $k$, $r$, and $\epsilon$.

By \eqnref{sy}, we have \be \label{dsfrac} \frac{\p u}{\p
\nu_B}\Big|_\pm = \frac{\p \Dcal_\Om f}{\p\nu_B}\Big|_\pm
-\frac{\p \Scal_\Om g}{\p\nu_B}\Big|_\pm +\frac{\p \Scal_B
\vp}{\p\nu_B}\Big|_\pm =(\lambda\pm\frac{1}{2})\vp
  \quad \mbox{on} \ \p B,
\ee and hence we obtain \be\label{x1est} \Big|  \frac{\p u}{\p
\nu_B}\Big|_+(X_1)\Big|  \geq\ds\frac{C\inf_{X\in
J_1}\Big|\frac{\p (\Dcal_\Om f)}{\p x}(X)\Big|}{1-\sigma +
4  r^*\sqrt{\ep}} .\ee

In order to prove \eqnref{est-100}, we now estimate from below
$g(X_2)$ where $X_2=(\rho,0)$, the point on $\p\Om$ which is the
closest to $\p B$. Elementary computations show that
\begin{align*}
&\frac{\p}{\p\nu_\Om} [(R_B R_\Om)^m\Dcal_\Om f](X_2) \\
&  = \prod^{m-1}_{n=0}g_\Omega(R_B (R_\Om R_B)^{n}(X_2))
g_B((R_\Omega R_B)^{n}(X_2)) \la \nabla \Dcal_\Om f((R_\Om
R_B)^m (X_2)), \nu_\Om(X_2) \ra .
\end{align*}
Since $(R_\Om R_B)^n (X_2)$ lies in $J_2$, the line between $X_2$
and $P_2$, for each $n$, we have as before
\begin{align*}
  \Big| g(X_2)\Big|  &\geq\Big|\frac{1}{2\lambda} \sum_{m=0}^{+\infty}
     \Big(\frac{a}{2\lambda}\Big)^m+2\Big|
     \cdot \inf_{X\in J_2}\Big|\frac{\p (\Dcal_\Om f)}{\p x}(X)\Big|
\\
   &\geq\ds\frac{\inf_{X\in J_2}\Big|\frac{\p (\Dcal_\Om f)}{\p
x}(X)\Big|}{1-\sigma +
   4 r^*\sqrt{\ep}}.
\end{align*}
Therefore, we get
 \be\label{x12est} \Big|  \frac{\p u}{\p
\nu_\Omega}\Big|_-(X_2)\Big|  \geq\ds\frac{C\inf_{X\in
J_2}\Big|\frac{\p (\Dcal_\Om f)}{\p x}(X)\Big|}{1-\sigma +
4 r^*\sqrt{\ep}} .\ee

  We now prove \eqnref{est1}. Let $N$ be the first
integer such that $N > \frac{1}{4 r^* \sqrt{\ep}}$. We then get
from Lemma \ref{lem304} that
\begin{align*}
&|\nabla ((R_\Om R_B)^m \Dcal_\Om f)(X)|  \\
&  \le \prod^m_{n=0} g_B(R_\Om(R_BR_\Om)^n(X))
g_\Om((R_BR_\Om)^n(X)) | \nabla \Dcal_\Om f((R_B R_\Om)^n(X))|  \\
& \le
\begin{cases}
\ds  \|\nabla \Dcal_\Om f \|_{L^\infty(\overline{\Om})} &\mbox{for
all } m, \\ \nm \ds  \|\nabla \Dcal_\Om f
\|_{L^\infty(\overline{\Om})} b^{m-N} \quad & \ds \mbox{if } m\geq
N,
\end{cases}
\end{align*}
where $b:= (1+r^* \sqrt{\epsilon})^{-1}$. It then follows that for
all $X\in \overline B$
\begin{align}
  \Big| \nabla\Scal_B\vp(X)\Big|  &\leq
\frac{1}{|\lambda|}\sum_{m=0}^{+\infty}
  \ds\frac{1}{(2|\lambda|)^m} \big |\nabla
  [(R_\Om R_B)^m\Dcal_\Om f](X)\big|\nonumber\\
  &\leq 2\|\nabla \Dcal_\Om f \|_{L^\infty(\Om)}\Big(\ds\sum_{m<
N }
     \Big(\frac{\ds 1}{\ds 2|\lambda|}\Big)^m + \frac{1}{|2\lambda|^N}
\sum_{m =0}^{+\infty}
     \Big(\frac{\ds b}{\ds 2|\lambda|}\Big)^m
    \Big)\nonumber \\
   &\leq C \|\nabla \Dcal_\Om f
\|_{L^\infty(\Om)}\Big(\frac{1-|\sigma|^{1/(
r^*\sqrt{\epsilon})}}{1-|\sigma|}+
   \frac{1}{1-|\sigma|
+ r^*\sqrt{\epsilon}}
    \Big)\nonumber\\
    &\leq\ds\frac{C \|\nabla \Dcal_\Om f
    \|_{L^\infty(\Om)}}{1-|\sigma|
+ r^*\sqrt{\epsilon}}\label{phiest}
    .
  \end{align}
By Lemma \ref{lm1} and \eqnref{rif}, $
|\nabla\Scal_B\vp(X)|=g_B(X)|\nabla\Scal_B\vp(R_B(X))|$ for $X\in
\Om\setminus \overline B$, and hence  \eqnref{phiest} holds for
all $X\in\Om$, {\em i.e}.,
 \be \label{phiest-2}
  \| \nabla\Scal_B\vp\|_{L^\infty(\Om)} \leq\ds\frac{C \|\nabla
 \Dcal_\Om f \|_{L^\infty(\Om)}}{1-|\sigma| + r^*\sqrt{\epsilon}} .
 \ee

Since $u$ is harmonic in $B$, it follows from \eqnref{dsfrac} that
 $$
 u(X)= - 2 \Scal_\Om \left (\frac{\p u}{\p \nu}\Big|_- \right )(X) + \mbox{constant}  (1-2\lambda) \Scal_B\vp (X) + \mbox{constant},
  \quad X \in  B.
 $$
Since $2\lambda-1= 2/(k-1)$,  \eqnref{phiest-2} gives
 \be\label{insidedisc}
 \| \nabla u \|_{L^\infty(B)} \le
 \ds\frac{C\|\nabla \Dcal_\Om f
    \|_{L^\infty(\Om)}}{|k-1|(1-|\sigma|
 + r^*\sqrt{\epsilon})}.
 \ee

By \eqnref{srgnew} and \eqnref{lm13},
 $$
 \Scal_\Om g = -2\sum_{m=1}^{+\infty} \frac{1}{(2\lambda)^m}
  (R_\Om R_B)^m
  \Dcal_\Om f - \Dcal_\Om f+ \mbox{constant}.
 $$
By Corollary \ref{cor45} and computations similar to those in
\eqnref{phiest}, we obtain \be\label{gest} \|\nabla \Scal_\Om g
\|_{L^\infty(\overline{\Om})} \le \ds\frac{C \|\nabla \Dcal_\Om f
\|_{L^\infty(\overline{\Om})}}{1-|\sigma| + r^*\sqrt{\epsilon}}.
\ee

By combining \eqnref{cd}, \eqnref{phiest}, and \eqnref{gest}, we
arrive at
 \be \label{last}
 \| \nabla u \|_{L^\infty(\Om)} \le C \frac{\| \nabla \Dcal_\Om f
 \|_{L^\infty(\Om)}}{1-|\sigma|
+ r^*\sqrt{\epsilon}}.
 \ee
Since $\Dcal_\Om$ maps $\mathcal{C}^{1,\alpha}(\p\Om)$ into
itself, it follows from the maximum principle that
 $$
 \|\nabla \Dcal_\Om f \|_{L^\infty(\overline{\Om})} \le \|\nabla
 \Dcal_\Om f \|_{L^\infty(\p\Om)} \le C\| f
 \|_{\mathcal{C}^{1,\alpha}(\p\Om)},
 $$
and hence we get \eqnref{est1}. The proof of Theorem
\ref{secondmain} is now complete. \qed

\noindent{\sl Proof of Theorem \ref{thirdmain}}. To prove Theorem
\ref{thirdmain} we use Lemma \ref{conj-2}. Let $G$ be the
$\mathcal{C}^{1,\alpha}$ function such that $\pd{G}{T}=g$ on
$\p\Om$ and $\int_{\p\Om} G =0$, and $v$ be the solution to
\eqnref{eqv1}. Since $\Dcal_\Om (G)$ is a harmonic conjugate to
$\Scal_\Om g$ in $\Om$ and $\RR^2 \setminus \overline{\Om}$ by
Lemma \ref{conj-2}, we have $$
 \la \nabla
 \Scal_\Om(g)(X), T_B(X_1) \ra= -\la \nabla \Dcal_\Om(G)(X),
 \nu_B (X_1) \ra,
 $$
 and
 $$
 \la \nabla
 \Scal_\Om(g)(X), T_\Omega(X_2) \ra= -\la \nabla \Dcal_\Om(G)(X),
 \nu_\Omega (X_2) \ra, \quad X \in \RR^2 \setminus \p\Om.
 $$
Thus \eqnref{est-neu} and \eqnref{est-neu-100} follow from
\eqnref{est}, \eqnref{est-100}, and \eqnref{onB}.

As one can see in the proof of Lemma \ref{conj-2}, $v$ is a
harmonic conjugate to $u$ in $\Om\setminus\overline B$ and
$\ds\frac{1}{k}v$ is a harmonic conjugate to $u$ in $B$. Therefore
by \eqnref{last} we get
 \be \|\nabla u
 \|_{L^\infty(\Om\setminus\overline{B})} \le \|\nabla v
 \|_{L^\infty(\Om)} \le \ds\frac{C \|\nabla \Scal_\Om g
 \|_{L^\infty(\Om)}}{1-|\sigma| + r^*\sqrt{\epsilon}}.
 \ee
Moreover,  by \eqnref{insidedisc}, we have
 \begin{align*}
 \|\nabla u
\|_{L^\infty(B)} & \le  \frac{1}{k}\|\nabla v \|_{L^\infty(B)} \\
&\le \frac{1}{k}\min\{\frac{1}{|\frac{1}{k}-1|},1\} \frac{C
\|\nabla \Scal_\Om g \|_{L^\infty(\Om)}}{1-|\sigma| +
r^*\sqrt{\epsilon}}\\&\le \frac{C \|\nabla \Scal_\Om g
\|_{L^\infty(\Om)}}{1-|\sigma| + r^*\sqrt{\epsilon}}.
 \end{align*}
Since $\|\nabla \Scal_\Om g \|_{L^\infty(\Om)} \le C \| g
\|_{\mathcal{C}^{\alpha}(\p\Om)}$, we finally get
\eqnref{est1-neu}. This completes the proof of Theorem
\ref{thirdmain}. \qed

\section*{Acknowledgments}  This work was partly supported by
CNRS-KOSEF grant No. 14889, ACI Nouvelles
    Interfaces des Math{\'e}matiques No. 171, and Korea Science and Engineering Foundation  grant
R02-2003-000-10012-0. Mikyoung Lim was partly supported by the
post-doctoral fellowship program of the Korea Science and
Engineering Foundation.

\end{document}